\documentclass[11pt,letterpaper]{article}
\usepackage[margin=1in]{geometry}
\usepackage{times}
\usepackage{amsmath,amssymb,amsthm,mathtools}
\usepackage{graphicx}
\usepackage{xcolor}
\usepackage{booktabs}
\usepackage{float}
\usepackage{natbib}
\usepackage[algo2e,ruled,vlined]{algorithm2e}
\usepackage{url}
\usepackage[
    colorlinks=true,
    linkcolor=blue,
    citecolor=blue,
    urlcolor=blue
]{hyperref}
\bibpunct{(}{)}{;}{a}{,}{,}
\SetAlgoVlined

\newtheorem{theorem}{Theorem}

\newtheorem{proposition}[theorem]{Proposition}

\newtheorem{corollary}[theorem]{Corollary}

\renewenvironment{proof}%
  {\par\noindent{\bfseries\upshape \proofname\ }}%
  {\hfill\rule{1.5ex}{1.5ex}\par\bigskip}

\DeclareMathOperator{\tr}{tr}

\newcommand{\ignore}[1]{}

\title{Mixed-Integer Programming for Causal Discovery with Unknown Hard-Intervention Targets}
\author{%
Tong Xu\thanks{Department of Industrial Engineering and Management Sciences, Northwestern University, Evanston, IL 60208, USA. Email: \texttt{tongxu2027@u.northwestern.edu}.}%
\qquad
Armeen Taeb\thanks{Department of Statistics, University of Washington, Seattle, WA 98195-4322, USA. Email: \texttt{ataeb@uw.edu}.}%
\qquad
Simge K\"u\c{c}\"ukyavuz\thanks{Department of Industrial Engineering and Management Sciences, Northwestern University, Evanston, IL 60208, USA. Email: \texttt{simge@northwestern.edu}.}%
\qquad
Ali Shojaie\thanks{Department of Biostatistics, University of Washington, Seattle, WA 98195-1617, USA. Email: \texttt{ashojaie@uw.edu}.}%
}
\date{}

\begin{document}

\maketitle

\begin{abstract}
When intervention targets are uncertain, learning causal relationships and
identifying directly perturbed variables become a joint task.
We develop an exact mixed-integer program based on optimized parent-set scores (PS-MIP) for joint graph and intervention-target estimation in
linear Gaussian structural equation models under stochastic hard interventions.
The method can incorporate candidate-edge restrictions and partial target
knowledge when available. For each random variable, PS-MIP evaluates candidate sets
of direct causes, or parent sets, by optimizing
the likelihood parameters and compatible target patterns, retaining both the
score and a minimizing target pattern. A global mixed-integer program then
selects parent sets subject to acyclicity. Under stated conditions, this
reduction preserves the joint objective and reconstructs the graph, targets,
and likelihood parameters. Because the local scores are precomputed, PS-MIP replaces a direct large bilinear mixed-integer program with a linear selection problem; empirically it certifies global optimality on far more instances than the direct formulation (126 of 180 versus 34) and is often orders of magnitude faster. On synthetic data, PS-MIP achieves favorable
false-discovery--true-discovery trade-offs relative to several benchmark algorithms for causal discovery with unknown intervention targets, along with strong intervention-target recovery in denser settings. An exploratory study
on real data from a controlled physical system further demonstrates recovery
beyond simulation.
\end{abstract}

\section{Introduction}

Learning causal relationships among variables is a central problem in fields
including computational biology, the social sciences, and economics
\citep{sachs2005causal,angrist1996identification,heckman2000causal}.  A
directed acyclic graph (DAG) represents direct causal relationships and
encodes conditional-independence restrictions on the observational
distribution. However, observational data alone may be compatible with
several DAGs that encode the same conditional-independence restrictions;
such DAGs belong to the same Markov equivalence class. Under the usual fully observed
Markov and faithfulness assumptions, and absent additional structural or
distributional restrictions, the population observational distribution
therefore generally identifies only this equivalence class rather than a
unique DAG \citep{verma1990equivalence}.

Interventions can reduce this ambiguity by changing how selected variables
are generated. We consider stochastic hard interventions: a targeted variable
receives an external random input in place of its dependence on its usual
direct causes, removing its incoming causal edges. The random-input
distribution may vary across environments
\citep{korb2004varieties,eberhardt2007interventions}. These interventions can substantially reduce the equivalence class  \citep{hauser2012characterization}.

In practice, the realized targets may not be fully known.  Gene-editing
experiments may exhibit off-target effects \citep{squires2020permutation},
while chemical perturbations may directly alter mechanisms beyond their
intended targets \citep{eaton2007exact}.  Since interventions are often unknown, we must jointly estimate the causal structure and the intervention targets.

We formulate this joint learning problem as an exact convex mixed-integer
program (MIP). The model selects a baseline DAG and an intervention-target
set for each environment while estimating shared and environment-specific
likelihood parameters. It readily incorporates prior candidate-edge restrictions
and partial target-status information. To simplify the global optimization,
we derive PS-MIP, an equivalent parent-set formulation based on precomputed
local scores. For each candidate parent set of a node, the continuous
likelihood parameters and compatible target patterns are optimized locally.
Precomputing these scores removes the likelihood parameters and target
indicators from the global optimization. The resulting MIP is a \emph{linear} mixed-integer program that selects parent sets that respect acyclicity. From the selected parents, we can reconstruct both the graph and the unknown interventional targets.

\begin{samepage}
\paragraph{Contributions.}
\begin{itemize}
    \item Under linear Gaussian structural causal models with unknown stochastic hard interventions, we give an exact MIP formulation for jointly learning the DAG and the intervention targets. The formulation allows for optional candidate-edge restrictions and partial target knowledge.
    \item We prove an exact reduction of the joint formulation to a parent-set selection \emph{linear} MIP (PS-MIP): the intervention targets and environment-specific likelihood parameters are eliminated in closed form into precomputed local scores, from which the graph, targets, and likelihood parameters are then reconstructed.
    \item We demonstrate graph and intervention-target recovery on synthetic
    data and graph recovery on real measurements from a controlled physical system.
    \item Empirically, this reduction is what makes exact solution practical: on the synthetic benchmarks PS-MIP certifies global optimality on 126 of 180 fits versus 34 for the direct joint formulation, and is up to orders of magnitude faster (e.g., 57\,seconds vs.\ the one-hour limit).
\end{itemize}
\end{samepage}

\subsection{Related work}

When intervention targets are known, intervention information can be incorporated through score-based or permutation-based search. Greedy Interventional Equivalence Search (GIES) extends greedy equivalence search to interventional data \citep{hauser2012characterization}, building on the characterization of interventional Markov equivalence classes from combined observational and interventional data \citep{hauser2015jointly}. Interventional Greedy Sparsest Permutation (IGSP) guides permutation search with interventional conditional-independence information \citep{wang2017permutation}, with the underlying equivalence classes characterized by \citet{yang2018characterizing}. When targets are unknown, Unknown-Target IGSP (UT-IGSP) jointly estimates the targets and the interventional equivalence class using conditional-independence and invariance information \citep{squires2020permutation}. Greedy Noise-Interventional Equivalence Search (GnIES) performs score-based equivalence search under Gaussian noise interventions with unknown targets \citep{gamella2022characterization}. Its intervention model retains the coefficients from a target's parents while allowing its noise distribution to change, in contrast to the stochastic hard interventions considered here, which remove incoming effects at targeted variables. Bayesian Causal Discovery with unknown interventions (BaCaDI) jointly approximates a posterior over graphs and intervention targets \citep{hagele2023bacadi}, and differentiable approaches have also been developed for interventional data with uncertain targets
\citep{ke2019learning,brouillard2020differentiable}.

Exact treatment of uncertain intervention targets predates this work.
\citet{eaton2007exact} extend exact Bayesian structure learning to account for uncertainty in intervention targets by summing over possible target assignments. Their approach relies on decomposable Bayesian models for which local marginal likelihoods can be computed tractably and retains the exponential subset enumeration underlying exact Bayesian dynamic programming. Under stochastic hard interventions, we instead use a penalized-likelihood estimator to select a graph and a target assignment, balancing fit to the data against penalties on the numbers of edges and intervention targets.

Exact integer optimization has an established role in observational Bayesian-network and Gaussian DAG learning, where mixed-integer formulations encode graph structure and acyclicity and provide globally valid bounds on the optimized objective
\cite[e.g.,][]{cussens2017bayesian,Manzour21,kucukyavuz2022consistent,xu2025integer}. For fully observed data, \cite{cussens2011bayesian} reduces structure learning to selecting one parent set per node from precomputed local scores. Related parent-set formulations also precompute local scores and impose acyclicity globally \citep{jaakkola2010learning}. This decomposition is specific to the observational likelihood. Under unknown interventions, the graph, the environment-specific mechanisms, and the latent target indicators enter jointly, so a direct encoding is a mixed-integer program with bilinear coefficient-switching terms (Appendix~\ref{app:reference-encoding}) rather than a parent-set score. Motivated by the use of local optimality information in mixed-integer optimization \citep{xu2026coordinate}, we solve the local fitting problems for each candidate parent set by eliminating the likelihood parameters in closed form and minimizing over compatible intervention-target patterns. This yields an exact parent-set formulation with precomputed scores despite the initially coupled graph, target, and parameter decisions.

A complementary literature exploits stability or structured changes across environments. Invariant causal prediction, BACKSHIFT, regression-invariance methods, and joint causal inference use different forms of invariance, distributional shift, or context augmentation to identify causal structure
\citep{peters2016causal,rothenhausler2015backshift,
ghassami2017learning,mooij2020joint}. A separate line of work reformulates structure learning as continuous, gradient-based optimization, enforcing acyclicity through a smooth algebraic constraint rather than a combinatorial one \citep{zheng2018dags,yu2019dag}, and this framework has been extended to interventional data with uncertain targets \citep{ke2019learning}, including Differentiable Causal Discovery from Interventional data (DCDI) \citep{brouillard2020differentiable}. These approaches optimize nonconvex continuous objectives and therefore typically yield only locally optimal solutions, whereas the formulation studied here targets the specified discrete score directly as a mixed-integer program, whose solver certifies global optimality, or reports a bounding optimality gap if stopped early. Prior work also shows that continuous methods can be sensitive to variance structure
\citep{reisach2021beware,seng2022tearing}.
\citet{xu2025asymptotically} highlight that similar objective values can correspond to substantially different Markov equivalence classes, and that asymptotic optimality alone does not establish statistical consistency. This motivates pursuing globally optimal solutions and using objective bounds to quantify any remaining optimization error. We report graph and intervention-target recovery alongside optimization diagnostics.

\section{Model and Joint Learning Objective}
\label{sec:int-formulation}
Consider \(p\) random variables, represented by the vertices
\(V=[p]:=\{1,\ldots,p\}\). A directed acyclic graph
(DAG) \(\mathcal G=(V,E)\) has directed edge set
\(E\subseteq\{(i,j)\in V\times V:i\neq j\}\), where \((i,j)\) represents
\(i\to j\). A centered random vector
\(X\in\mathbb R^p\) follows the baseline linear Gaussian structural equation
model
\[
    X=B^\top X+\epsilon,
    \qquad
    \epsilon\sim\mathcal N(0,\Omega),
\]
where \(B\in\mathbb R^{p\times p}\) satisfies \(B_{ij}=0\) whenever
\((i,j)\notin E\), and \(\Omega\) is positive definite and diagonal.  Thus,
column \(j\) of \(B\) contains the incoming coefficients of node \(j\). The matrix $\Omega$ being diagonal ensures that there is is no latent confounding (i.e., causal sufficiency).

Let \(\mathcal E\) be the finite environment set, where \(e_0
\in \mathcal E\) represents the observational
environment, and \(\mathcal E_{\setminus0}=\mathcal E\setminus\{e_0\}\) represents environments that interventions have occured.
Each environment $e$ has a target set \(I_e\subseteq[p]\), with
\(I_{e_0}=\varnothing\). A stochastic hard intervention replaces each
targeted equation by independent centered Gaussian noise of variance
\(\omega_j^e>0\), deleting its incoming coefficients
\citep{korb2004varieties}. Non-target nodes retain both baseline coefficients
and noise variances. 

Let \(\Sigma^e\) denote the covariance in environment \(e\).
We parameterize its inverse as \((\Sigma^e)^{-1}=\Gamma^e{\Gamma^e}^{\top}\),
where \(\Gamma^e\) is a \emph{precision factor} with positive diagonal.
In the observational environment, write
\(\Gamma^0:=\Gamma^{e_0}=(\mathbf I_p-B)\Omega^{-1/2}\), where \(\mathbf I_p\) is the
\(p\times p\) identity matrix. Column \(j\) of \(\Gamma^0\) contains the baseline residual
coefficients divided by the noise standard deviation:
\(\Gamma^0_{jj}=1/\sqrt{\Omega_{jj}}\) and
\(\Gamma^0_{ij}=-B_{ij}/\sqrt{\Omega_{jj}}\) for \(i\neq j\).
Thus \(\Gamma^0\) and \(B\) have identical off-diagonal zero patterns;
their nonzero off-diagonal entries encode the same directed edges.
Across environments, the sharing rule is
\begin{equation}
    \Gamma^e_{\cdot j}
    =
    \begin{cases}
        \Gamma^0_{\cdot j}, & j\notin I_e,\\
        (\omega_j^e)^{-1/2}\mathbf e_j, & j\in I_e,
    \end{cases}
    \qquad e\in\mathcal E,\ j\in[p].
    \label{eq:mechanism-sharing}
\end{equation}
Here \(\mathbf e_j\) is the \(j\)th standard basis vector. For a targeted
node \(j\), \(\omega_j^e\) is the variance of the independent Gaussian input
that replaces its usual mechanism in environment \(e\).
Appendix~\ref{app:mechanism-representation} expands the equivalent
coefficient, variance, and covariance representation.

For each \(e\in\mathcal E\), let
\(\mathcal X^e\in\mathbb R^{n_e\times p}\) contain \(n_e\) independent
centered observations.  Define
\(\widehat\Sigma^e={{\mathcal X^e}^{\top}\mathcal X^e}/{n_e}\), let
\(n=\sum_{h\in\mathcal E}n_h\), and set \(w_e=n_e/n\).
For notational simplicity, we represent the collection $\{\Gamma^e\}_{e\in\mathcal E}$ by $\Gamma$.
Up to an additive constant and a common positive factor, the joint
normalized negative log-likelihood is 
\begin{equation}\label{eq:Lgamma}
    \mathcal L(\Gamma)
    :=\sum_{e\in\mathcal E}w_e
    \left\{
        \sum_{j=1}^p-2\log\Gamma^e_{jj}
        +\operatorname{tr}
        \bigl(\Gamma^e{\Gamma^e}^{\top}\widehat\Sigma^e\bigr)
    \right\}.
\end{equation}

The target sets may be known, partially known, or unknown.  For
\(e\in\mathcal E_{\setminus0}\), introduce
\(\delta_j^e=\mathbb I[j\in I_e]\) for \(j\in[p]\), where \(\mathbb I[\cdot]\) is the indicator function, equal to 1 when its argument holds and 0 otherwise.
Known target-status information fixes the corresponding indicators.
In the fully unknown setting, all are free.

Let \(E_{\mathrm{super}}\) denote the candidate directed edge set.
By default, \(E_{\mathrm{super}}=\{(i,j)\in V\times V:i\neq j\}\).
Prior structural information or a screening rule can restrict this set,
giving a candidate superstructure. Fix an edge
penalty \(\lambda_g\ge0\), and a target penalty \(\lambda_\delta\ge0\).
Let \(\mathfrak M\) be the search space of
admissible graphs, intervention-target assignments, and likelihood parameters.
It consists of triples \((\mathcal G,I,\Gamma)\) satisfying the following requirements:
\(\mathcal G=([p],E)\) is a DAG with \(E\subseteq E_{\mathrm{super}}\);
\(I=(I_e)_{e\in\mathcal E}\) respects the fixed target-status information and
\(I_{e_0}=\varnothing\); the columns of \(\Gamma\) obey the mechanism-sharing
rule above; and
\begin{equation*}
\begin{gathered}
    \Gamma^{e_0}_{ij}=0\quad(i\neq j,\ (i,j)\notin E),\\
    0<\Gamma^e_{jj}
    \quad(e\in\mathcal E,\ j\in[p]).
\end{gathered}
\end{equation*}
The edge set specifies allowed baseline entries, which may have zero fitted
values. The bounds apply to \(\Gamma\), with \(B\) recovered from its
noise-scaled columns. Diagonal entries must be positive for the logarithm
to be defined; the negative-log term diverges as an entry approaches zero.

The joint learning objective is
\begin{equation}
    \min_{(\mathcal G,I,\Gamma)\in\mathfrak M}
    \mathcal L(\Gamma)+\lambda_g|E|
    +\lambda_\delta\sum_{e\in\mathcal E_{\setminus0}}|I_e|.
    \label{eq:joint-learning-objective}
\end{equation}
Here \(|E|\) is the number of selected edges in \(\mathcal G\) and \(\sum_{e\in\mathcal E_{\setminus0}}|I_e|\) is the total number of intervention targets across environments; \(\lambda_g\ge0\) and \(\lambda_\delta\ge0\) are regularization parameters that promote graph sparsity and intervention sparsity, respectively.

Building on the observational mixed-integer formulation of
\citet{xu2025integer}, we obtain a direct formulation of the joint learning
problem by introducing intervention-target indicators and enforcing
mechanism sharing across environments. Appendix~\ref{app:reference-encoding}
gives its full mixed-integer convex encoding, Problem~\eqref{Problem:interventional},
using edge indicators, order variables, and bounded switching relations. Its
projection is precisely \(\mathfrak M\), with the same objective and bounds.
The direct formulation retains continuous fitting variables across environments 
alongside graph and target indicators. As a result, this formulation is slow to solve and does not scale to moderate problem sizes. The reformulation described next, PS-MIP, simplifies the global problem by replacing these fitting problems
with precomputed scores. This results in a much faster procedure for our experiments; Section~\ref{sec:synthetic-experiments} and Table~\ref{tab:joint-mip-comparison} quantify this gap.

\section{PS-MIP: Learning with Optimized Parent-Set Scores}
\label{sec:precomputed-parent-set}

Unlike the observational case, the unknown-target likelihood does not decompose into parent-set scores directly: each node's contribution also depends on its unknown intervention targets, which are coupled with the rest of the estimation. We develop a novel reformulation that restores this decomposition by eliminating the intervention targets and environment-specific parameters in closed form for each candidate parent set, so that the joint problem reduces exactly to a parent-set selection MIP with precomputed local scores. PS-MIP exploits the nodewise structure of the DAG likelihood. Once a node's candidate parents and its target pattern are specified, its likelihood contribution can be fitted independently of other nodes' parameters, using measurements of the node and
its parents across environments. We therefore optimize the likelihood
parameters and compatible target patterns for each candidate parent set,
retaining its score and a minimizing pattern. A global MIP selects one parent
set per node subject to acyclicity. Under the stated conditions, this
construction preserves the joint objective and allows reconstruction of the
graph, targets, and likelihood parameters.

\subsection{Optimized local scores}
\label{sec:local-scores}

For a proposed parent set of node \(j\), we seek the parameters and
intervention targets that best explain its observations across environments,
balancing data fit against the penalties in the joint objective. We do this
in two stages: fit the
parameters for a fixed target pattern, then compare the compatible patterns.
The parent set stays fixed throughout this calculation; parent sets are
selected globally in the next subsection.

For node $j$, let \(\operatorname{pa}^{\mathrm{sup}}(j)=\{i:(i,j)\in E_{\mathrm{super}}\}\).
Fix a candidate parent set \(A\subseteq\operatorname{pa}^{\mathrm{sup}}(j)\),
so baseline column entries outside \(A\cup\{j\}\) are zero. Temporarily also
fix \(S\subseteq\mathcal E_{\setminus0}\), the environments where node \(j\)
is targeted. This divides the fitting problem into two parts. In
\(N(S):=\{e_0\}\cup(\mathcal E_{\setminus0}\setminus S)\), the node follows
one shared regression on \(A\), with common coefficients and noise variance.
In \(S\), hard interventions remove its dependence on its parents, leaving
a separate noise-variance fit for each environment.

\paragraph{Columnwise likelihood.}
Both parts can be fitted independently of other nodes' parameters because the
likelihood \(\mathcal L(\Gamma)\) in ~\eqref{eq:Lgamma} separates
by node. Its log terms are already columnwise, and the trace term satisfies
\(\operatorname{tr}\!\big(\Gamma^e\Gamma^{e\top}\widehat\Sigma^e\big)
=\sum_{j=1}^p(\Gamma^e_{\cdot j})^\top\widehat\Sigma^e\,\Gamma^e_{\cdot j}\).
Collecting the terms in column \(j\), environment \(e\) contributes the
\emph{columnwise likelihood}
\begin{equation}\label{eq:columnwise-likelihood}
  \mathcal L^e_j(\Gamma^e_{\cdot j})
  = w_e\Big\{-2\log\Gamma^e_{jj}
     + (\Gamma^e_{\cdot j})^\top\widehat\Sigma^e\,\Gamma^e_{\cdot j}\Big\}.
\end{equation}
We minimize the total column-\(j\) cost over the two groups of environments
separately.

\paragraph{Fitting the shared regression.}
In the environments \(N(S)=\{e_0\}\cup(\mathcal E_{\setminus0}\setminus S)\),
node \(j\) is not intervened, so it shares one baseline column
\(\Gamma^0_{\cdot j}\) supported on \(A\cup\{j\}\). Summing
\eqref{eq:columnwise-likelihood} over these environments and using
\(\Gamma^0_{ij}=0\) for \(i\notin A\cup\{j\}\), the pooled cost depends on the
data only through the total weight and the pooled second-moment matrix
\begin{equation*}
  a_S=\sum_{e\in N(S)}w_e,
  \qquad
  Q_{jAS}=\sum_{e\in N(S)}w_e\,\widehat\Sigma^e_{(A,j),(A,j)}
        =\begin{pmatrix}Q_{AA}&Q_{Aj}\\ Q_{jA}&Q_{jj}\end{pmatrix},
\end{equation*}
where \((A,j)\) orders the parent coordinates first and node \(j\) last. Writing
\(t=\Gamma^0_{jj}>0\) for the inverse noise standard deviation and
\(u=\Gamma^0_{Aj}\) for the noise-scaled incoming coefficients, the pooled
column-\(j\) objective is
\begin{equation*}
  \Phi_{jAS}(u,t)
  =-2a_S\log t+\begin{pmatrix}u\\ t\end{pmatrix}^{\!\top}Q_{jAS}
     \begin{pmatrix}u\\ t\end{pmatrix}
  =-2a_S\log t+u^\top Q_{AA}u+2t\,u^\top Q_{Aj}+Q_{jj}t^2 .
\end{equation*}

\emph{Step 1 (coefficients).} Fix \(t\) and minimize \(\Phi_{jAS}\) over \(u\).
For \(A\neq\varnothing\) with \(Q_{AA}\) nonsingular, the first-order condition
\(Q_{AA}u+t\,Q_{Aj}=0\) gives \(u^\star(t)=-Q_{AA}^{-1}Q_{Aj}\,t\). Substituting
reduces \(\Phi_{jAS}\) to the scalar \(-2a_S\log t+r_{jAS}t^2\), where
\begin{equation*}
  r_{jAS}=Q_{jj}-Q_{jA}Q_{AA}^{-1}Q_{Aj}
\end{equation*}
is the Schur complement of \(Q_{AA}\) in \(Q_{jAS}\) (for \(A=\varnothing\), set
\(r_{jAS}=Q_{jj}\)). Here \(r_{jAS}/a_S\) is the pooled mean-squared residual of
node \(j\) given \(A\); we assume \(r_{jAS}>0\), for which \(Q_{jAS}\succ0\) is
sufficient.

\emph{Step 2 (noise scale).} Minimizing the reduced objective over \(t>0\), the
condition \(-2a_S/t+2r_{jAS}t=0\) gives the baseline optimum
\begin{equation*}
  \Gamma^{0,\star}_{jj}=\sqrt{\tfrac{a_S}{r_{jAS}}},
\end{equation*}
and back-substitution into Step~1 gives the coefficients
\begin{equation*}
  \Gamma^{0,\star}_{Aj}=-Q_{AA}^{-1}Q_{Aj}\,\Gamma^{0,\star}_{jj}.
\end{equation*}
Evaluating \(\Phi_{jAS}\) at \((u^\star,t^\star)\) yields the minimized
shared-mechanism cost
\begin{equation}\label{eq:shared-cost}
  s^{\mathrm{sh}}_{jAS}=a_S\Big[\,1+\log\!\big(r_{jAS}/a_S\big)\Big].
\end{equation}

\paragraph{Fitting the targeted environments.}
For each \(e\in S\) the hard intervention deletes node \(j\)'s incoming edges, so
its column reduces to \(\Gamma^e_{\cdot j}=t_e\,\mathbf e_j\) with only the noise
scale \(t_e=\Gamma^e_{jj}>0\) left to fit. Applying the columnwise likelihood
\eqref{eq:columnwise-likelihood} to this column,
\begin{equation*}
  \mathcal L^e_j=w_e\big\{-2\log t_e+\widehat\Sigma^e_{jj}\,t_e^2\big\}.
\end{equation*}
Minimizing over \(t_e>0\) gives the intervened optimum
\(\Gamma^{e,\star}_{jj}=1/\sqrt{\widehat\Sigma^e_{jj}}\) with cost
\begin{equation}\label{eq:targeted-cost}
  s^{\mathrm{tar}}_{ej}=w_e\big(1+\log\widehat\Sigma^e_{jj}\big),
\end{equation}
which requires \(\widehat\Sigma^e_{jj}>0\). These fits are independent across
\(e\in S\).

\paragraph{Local score.}
Combining the shared cost \eqref{eq:shared-cost} with the targeted costs
\eqref{eq:targeted-cost}, and charging \(\lambda_g\) for each selected parent and
\(\lambda_\delta\) for each targeted environment, the local score for the triple
\((j,A,S)\) is
\begin{equation*}
\begin{aligned}
  s(j,A,S)
  &= s^{\mathrm{sh}}_{jAS}+\lambda_g|A|
     +\sum_{e\in S}\big(s^{\mathrm{tar}}_{ej}+\lambda_\delta\big)\\
  &= a_S\Big[1+\log\!\big(r_{jAS}/a_S\big)\Big]+\lambda_g|A|
     +\sum_{e\in S}\Big[w_e\big(1+\log\widehat\Sigma^e_{jj}\big)+\lambda_\delta\Big].
\end{aligned}
\end{equation*}

The first term is the shared-regression cost, \(\lambda_g|A|\) penalizes
edges, and the sum fits and penalizes targeted environments. Changing the
target pattern therefore changes both which observations are pooled and
the cost of explaining the remaining environments separately.

\paragraph{Choosing compatible target patterns.}
With the continuous parameters fitted, we can now compare target patterns
while keeping the parent set \(A\) fixed.
Let \(\mathcal S_j\) contain the patterns consistent with any fixed target
status of node \(j\). With fully unknown targets,
\(\mathcal S_j\) is the power set of \(\mathcal E_{\setminus0}\).
We call target choices \emph{nodewise feasible} when each node's allowed
pattern can be chosen separately: any collection of choices from the
\(\mathcal S_j\) is jointly feasible. Formally, the feasible collection is
\(\prod_j\mathcal S_j\). Fixed individual target statuses satisfy this
condition, since they restrict only the corresponding node's patterns.
Under this condition, the lowest-cost pattern for \(A\) can be chosen locally:
\begin{align*}
    s(j,A)&=\min_{S\in\mathcal S_j}s(j,A,S),
    \\
    S^*(j,A)&\in\operatorname*{arg\,min}_{S\in\mathcal S_j}s(j,A,S).
\end{align*}
Retain any minimizing pattern with each score. Repeating this calculation
over candidate parent sets gives the inputs to the global MIP: one optimized
score and one target pattern for each choice. Different parent sets can
retain different target patterns, so selecting parent sets in the next
subsection determines both the graph and the intervention-target estimates.

\subsection{Global selection, equivalence, and reconstruction}

Choosing each node's lowest-cost parent set independently can produce a
directed cycle. We therefore select parent sets jointly, minimizing their
total precomputed cost subject to acyclicity.

Let \(x_{jA}=1\) if parent set \(A\) is selected for node \(j\), and let
\(x_{jA}=0\) otherwise. For each node, a continuous ordering variable
\(\psi_j\) is constrained to increase along selected edges. The MIP is
\begin{subequations}
\begin{align}
    \min_{x,\psi}\quad
    &\sum_{j=1}^{p}
      \sum_{A\subseteq\operatorname{pa}^{\mathrm{sup}}(j)}
      s(j,A)x_{jA}
    \label{eq:parent-set-mip-objective}\\
    \text{s.t.}\quad
    &\sum_{A\subseteq\operatorname{pa}^{\mathrm{sup}}(j)}
      x_{jA}=1,
      &&j=1,\ldots,p,
    \label{eq:one-parent-set}\\
    &1-p
      +p\sum_{\substack{
          A\subseteq\operatorname{pa}^{\mathrm{sup}}(j)\\
          i\in A}}
      x_{jA}
      \leq \psi_j-\psi_i,
      &&(i,j)\in E_{\mathrm{super}},
    \label{eq:parent-set-order}\\
    &1\leq\psi_j\leq p,
      &&j=1,\ldots,p,
    \label{eq:parent-set-order-bounds}\\
    &x_{jA}\in\{0,1\},
      &&j=1,\ldots,p,\quad
        A\subseteq\operatorname{pa}^{\mathrm{sup}}(j).
    \label{eq:parent-set-binary}
\end{align}
\label{eqn:MIP_form}
\end{subequations}
Constraint~\eqref{eq:one-parent-set} selects exactly one parent set for each
node. If the selected parent set for node \(j\) contains \(i\),
constraint~\eqref{eq:parent-set-order} reduces to
\(\psi_j\geq\psi_i+1\); hence the selected edges cannot form a directed
cycle. Separate edge variables are unnecessary because the term
\(\lambda_g|A|\) is already included in every local score.

\begin{proposition}[Equivalence of the parent-set formulation]
\label{prop:parent-set-equivalence}
Assume the local-score conditions in
Appendix~\ref{app:local-score-conditions} and that
feasible target patterns form \(\prod_j\mathcal S_j\). The parent-set MIP \eqref{eqn:MIP_form} has
the same optimal value as the joint learning problem
\eqref{eq:joint-learning-objective}, equivalently its reference encoding
\eqref{Problem:interventional}. An optimal parent-set solution reconstructs
an optimal graph, its intervention targets, and the corresponding likelihood
parameters within the same candidate class and bounds.
\end{proposition}
The proof is given in
Appendix~\ref{app:parent-set-equivalence-proof}.

For selected parent sets \(\widehat A_j\), recover
\(\widehat I_e=\{j:e\in S^*(j,\widehat A_j)\}\) and the associated local parameters.
The reconstruction in
Appendix~\ref{app:parent-set-equivalence-proof} gives the full mapping.
The correspondence uses the same candidate class and \(\Gamma\) bounds.
Thus, selecting the parent sets determines the graph, intervention targets,
and fitted parameters without further joint optimization. The recovered DAG and estimated intervention targets can also be converted into an interventional Markov equivalence class of DAGs \citep{hauser2012characterization}.

\subsection{Computation and screening}

With \(K=|\mathcal E_{\setminus0}|\) and
\(d_j=|\operatorname{pa}^{\mathrm{sup}}(j)|\), enumeration considers at most
\(\sum_j2^{d_j+K}\) parent-set/target-pattern combinations. The MIP \eqref{eqn:MIP_form}
 has
\(\sum_j2^{d_j}\) binary indicators and \(p\) continuous order variables.
Although enumeration is exponential in candidate parents and environments,
screening can substantially reduce the candidate collections for larger
graphs by reducing the \(d_j\). Our synthetic implementation uses an
observational graphical-lasso screen and admits both directions of each
retained pair. Partial target information further reduces the patterns
enumerated for each node.
Appendix~\ref{app:screening-justification} gives the screening argument, and
Appendix~\ref{app:parent-set-implementation} provides the full algorithm
and implementation details.

\section{Empirical Evaluation}
\label{sec:experiments}

\subsection{Synthetic experiments}
\label{sec:synthetic-experiments}

\paragraph{Setup.}
For each combination of \(p\in\{10,20,30\}\) variables and edge-density parameter \(e\in\{1,2\}\), we generate ten independent random DAGs and sampling edges, with \(\mathbb E[|E|]=ep\). Each of the 60 instances contains
1000 observational samples and five hard-intervention environments with
200 samples each. We compare unknown-target PS-MIP, UT-IGSP
\citep{squires2020permutation}, and GnIES \citep{gamella2022characterization}
with oracle PS-MIP, IGSP, and GIES \citep{hauser2012characterization},
which receive the generating targets. DCDI-G, the Gaussian-likelihood variant of DCDI \citep{brouillard2020differentiable}, including its oracle condition, and BaCaDI \citep{hagele2023bacadi} returned no finite accuracy metrics under the one-hour fitting budget per setting, and therefore do not appear in the plots
PS-MIP uses observational screening and sets the graph and target penalties equal, \(\lambda_g=\lambda_\delta\), at each tested setting.
Appendix~\ref{app:experimental-protocol} gives the complete data-generation,
parameter setting, and computational protocol. Appendix Table~\ref{tab:synthetic-runtime} reports median fitting times
by graph size and density. PS-MIP's median fitting time remains below
one minute in the sparser settings (\(e=1\)) but exceeds one hour in
the larger, denser settings (\(p\in\{20,30\}\), \(e=2\)).
To isolate the benefit of the reduction itself, we solve the same 60 instances and penalties with the direct joint MIP of Appendix~\ref{app:reference-encoding} under an identical budget. PS-MIP reaches certified optimality---a solver-proven guarantee that no better solution exists---on 126 of 180 fits against 34 for the joint formulation, and is markedly faster where both finish---1.43\,s vs.\ 577\,s at $(p,e)=(10,1)$, and 57\,s vs.\ the one-hour limit at $(30,1)$ (Table~\ref{tab:joint-mip-comparison}). Precomputing local scores is what makes exact joint estimation tractable at these sizes.

\begin{figure}[t]
    \centering
    \includegraphics[width=0.75\linewidth]{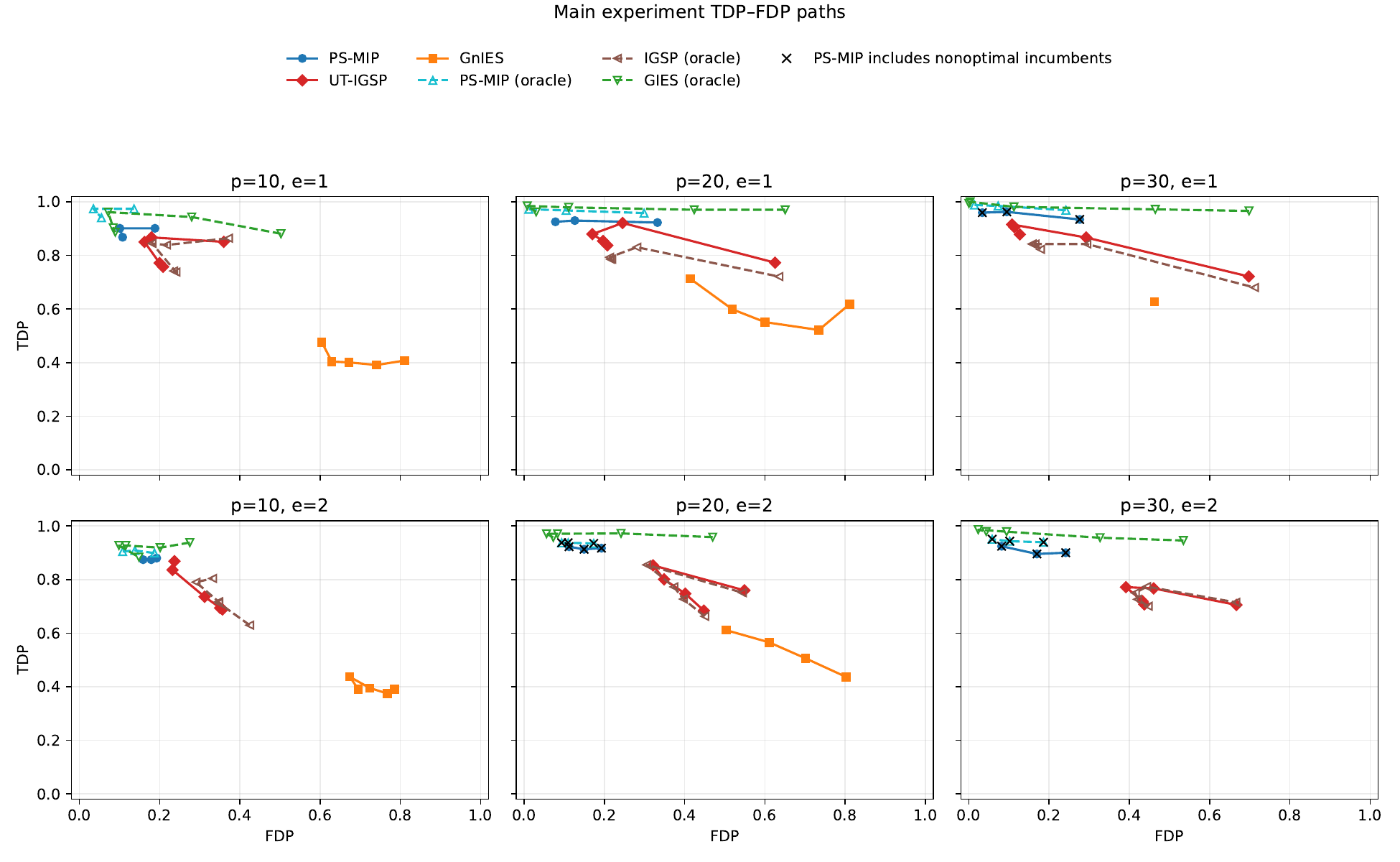}
    \caption{Synthetic graph-recovery FDP--TDP paths; \(p\) is the number of
    variables and \(e\) sets the expected edge count to \(ep\).
    Each point averages available replicates for its method
    and setting. Crosses mark PS-MIP settings containing at least one
    nonoptimal incumbent.}
    \label{fig:synthetic-tdp-fdp}
\end{figure}

\paragraph{Graph recovery.}
We report performance across a grid of regularization for the comparison methods. Detailed chosen grid for each method is given in Appendix~\ref{app:experimental-protocol}. Figure~\ref{fig:synthetic-tdp-fdp} plots the true-discovery proportion (TDP)
against the false-discovery proportion (FDP), evaluated over the interventional Markov equivalence class. TDP measures the fraction of true edges recovered, while FDP measures
the fraction of estimated edges that are false discoveries. Both use
worst-case comparisons across compatible DAG representations, with formal
definitions in Appendix~\ref{app:graph-metrics}. Curves nearer the upper-left corner combine
more true-edge recovery with fewer false discoveries.

PS-MIP shows favorable trade-offs relative to UT-IGSP, most clearly for
\(p\in\{20,30\}\). GnIES's returned paths have lower TDP and higher FDP,
consistent with a mismatch in intervention assumptions: its noise-intervention
model preserves parental effects, whereas our hard interventions remove them
\citep{gamella2022characterization}. Oracle GIES and oracle PS-MIP lie close to the upper-left corner
across most cells. Together, these comparisons show effective graph recovery
with unknown targets and illustrate the value of supplied target information.

\begin{figure}[t]
    \centering
    \includegraphics[width=0.75\linewidth]{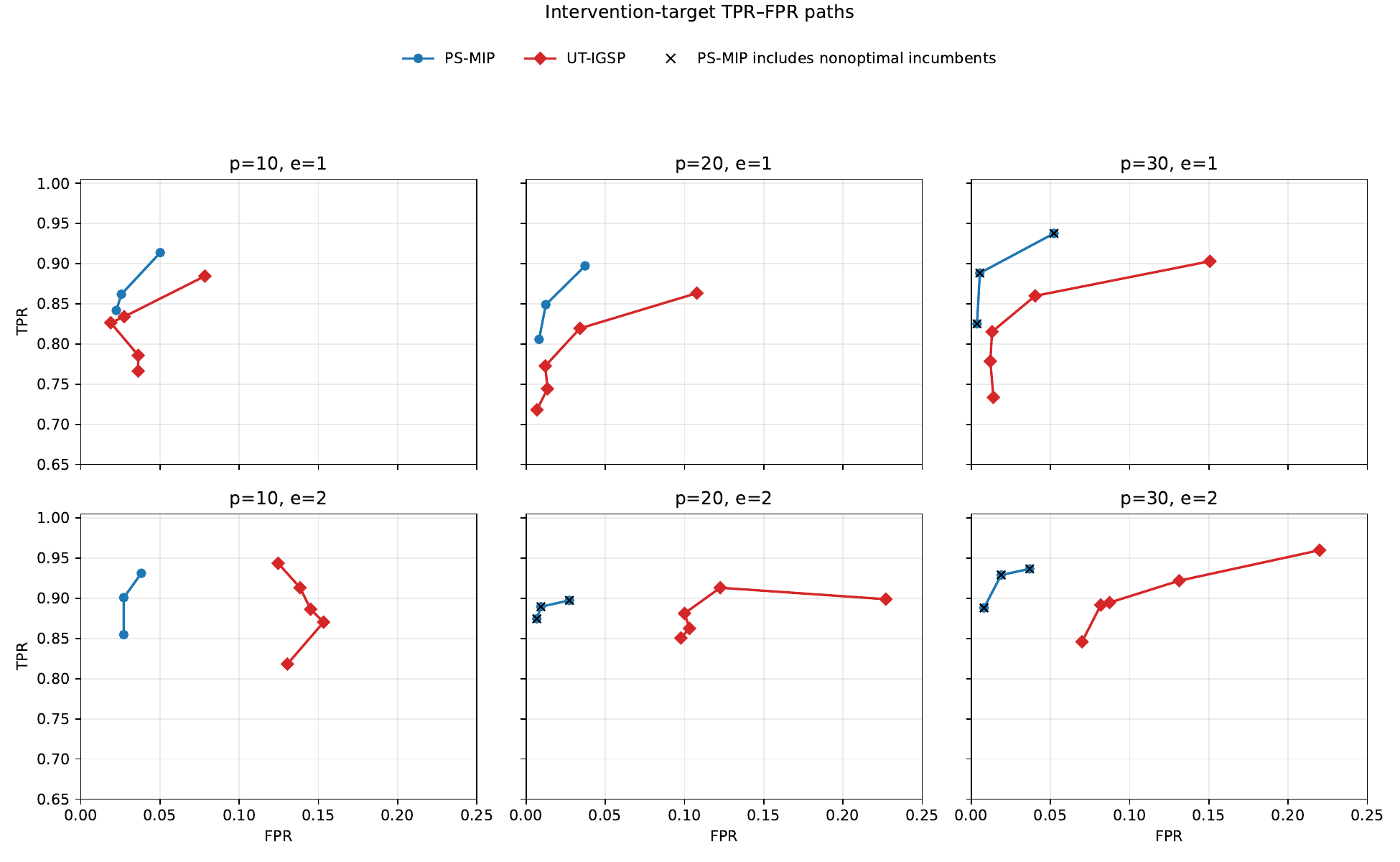}
    \caption{Synthetic intervention-target TPR--FPR paths. Crosses mark PS-MIP settings containing at least
    one nonoptimal incumbent.}
    \label{fig:target-tpr-fpr}
\end{figure}

\paragraph{Intervention-target recovery.}
Figure~\ref{fig:target-tpr-fpr} compares how accurately PS-MIP and UT-IGSP
identify the variables targeted in each environment. The true-positive rate
(TPR) measures the fraction of targets correctly identified, while the
false-positive rate (FPR) measures the fraction of non-targets incorrectly
classified as targets. Curves closer to the upper-left corner therefore
indicate better recovery. PS-MIP compares favorably with UT-IGSP across all
six settings. For denser graphs (\(e=2\)), PS-MIP has lower FPR across the
reported paths, while the TPR ranges overlap. For sparser graphs (\(e=1\)),
the methods perform more similarly, with PS-MIP generally recovering more
targets at similar or lower FPR. These findings complement the graph-recovery
results. 

\subsection{Real-data experiment}
\label{sec:chambers}

We use the Causal Chambers light tunnel, a physical setup with an
adjustable red--green--blue (R--G--B) LED source and sensors that measure its
electrical current and light output \citep{gamella2024causal}. Our analysis
uses configuration \texttt{scm\_4} of \texttt{lt\_camera\_v1}. Each row records
one measurement after the LED settings are applied. Its ten entries are
the three brightness settings \(R,G,B\), one current-sensor reading,
and six light-intensity readings: visible and infrared channels at three
sensor locations. The data combine the chosen brightness settings
with measurements of how the system responds.

In the reference environment, the brightness settings follow a programmed
chain \(R\to G\to B\): green depends on red, and blue depends on green.
Three other environments intervene on red, green, or blue by setting the
chosen channel independently of its usual inputs, using a new distribution.
For example, intervening on green removes its dependence on red, while blue
continues to depend on green and the sensors record how the system responds.
Each intervention cuts the causal links into its target.
We use all 40,000 observations, with 10,000 from each
environment. The reference graph combines the two programmed links with
21 documented physical links from the LED settings to the sensor readings.
To test learning with unknown targets, we give PS-MIP, UT-IGSP, and GnIES
the environment groups and identify the observational environment, but do
not tell them which channel was intervened on. PS-MIP considers all possible
edges. All methods use the same scaling based on the reference environment,
as described in Appendix~\ref{app:chambers}.

Figure~\ref{fig:chamber-unknown} compares true positives (TP) and false positives (FP) against
the reference graph. The directed panel evaluates saved DAG arcs, while
the skeleton panel evaluates adjacencies regardless of direction.
PS-MIP recovers all 23 reference directions at every tested graph penalty.
At its sparsest point, both panels give 23 TP and 4 FP, recovering all
21 physical and both programmed links. The smallest directed FP counts are
seven for UT-IGSP and six for GnIES. This result illustrates sparse
reference-graph recovery from real measurements with unknown targets.

\begin{figure}[t]
    \centering
    \includegraphics[width=0.65\linewidth]{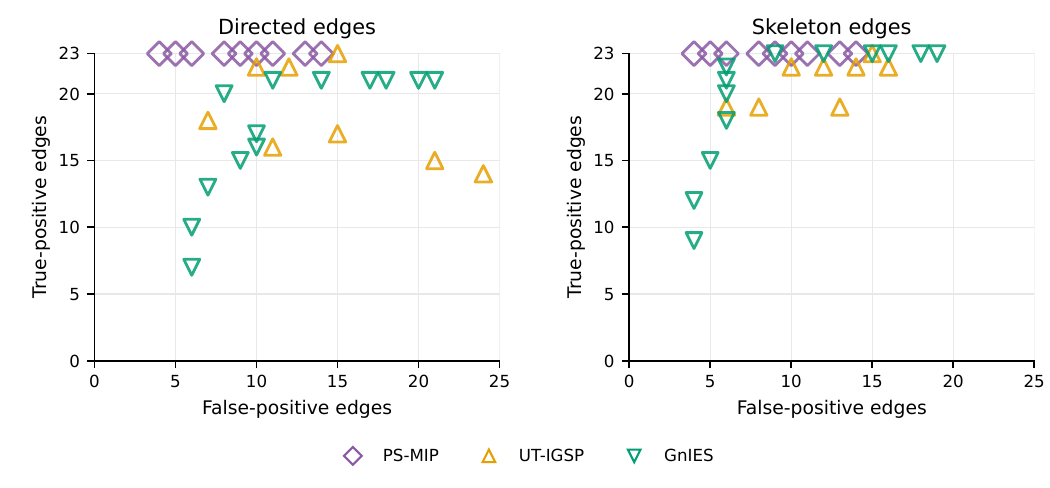}
    \caption{Unknown-target recovery on Causal Chambers \texttt{scm\_4}.
    Directed and skeleton TP--FP counts use the 23-edge reference DAG.
    Exact duplicate coordinates are drawn once within each method and panel.}
    \label{fig:chamber-unknown}
\end{figure}

\section{Discussion}
\label{sec:discussion-limitations}

In this paper, we bring causal graph learning and intervention-target estimation into
a single optimization problem. Under the conditions of Proposition~\ref{prop:parent-set-equivalence}, this representation preserves the joint objective. The same formulation accommodates unknown targets, partial target knowledge, and optional
restrictions on candidate edges.

The experiments demonstrate the value of this joint approach. On the synthetic
benchmarks, PS-MIP shows generally favorable graph- and target-recovery
trade-offs relative to UT-IGSP. The Causal Chambers study further demonstrates
sparse reference-graph recovery from real measurements with unknown targets.
Together, these findings support the approach on the moderate-size problems
studied here. The evidence spans a linear Gaussian hard-intervention design
and one exploratory physical configuration, motivating evaluation in a wider
range of settings.

Three priorities for future work are scalability, broader intervention
models, and extending beyond linear models. Local enumeration grows exponentially with the number of candidate
parents and environments, making screening that preserves true edges and
more efficient score computation important next steps. Extending the local
scores to interventions that change, rather than remove, the effects of a
variable's parents would broaden the model's scope. Building on the MIP formulation in \cite{Zhang2025ConvexMP}, it would be of interest to extend our approach to non-linear settings. Evaluation across more
intervention regimes and physical systems would help assess these extensions.
Overall, this work establishes optimized parent-set scores as a concrete
approach to joint causal graph and intervention-target learning.

\subsection*{AI use statement}
We used OpenAI Codex to assist with code, plotting, and \LaTeX{} editing. Physical measurements come from
the cited public dataset, and model fits are produced by the stated algorithms and solver.  All AI-assisted work is subject to author review; the authors remain responsible for the experiments, mathematical claims, and manuscript.

\subsection*{Ethics statement}
This work uses synthetic data and public physical-sensor data from the
Causal Chambers, with attribution to their creators. It involves no new
human-participant recruitment, physical data collection, or collection of
private data.
Consequential uses should include domain-expert validation and assessment of
harms from incorrect causal conclusions.

\subsection*{Reproducibility statement}
Sections~\ref{sec:int-formulation}--\ref{sec:experiments}
specify the joint objective, PS-MIP, data generation, and evaluation. Appendices
\ref{app:synthetic-diagnostics} and \ref{app:chambers}
give settings and diagnostics. The Chamber pipeline preserves all 63 final
fits, target-information conditions, common-input hashes, source revisions,
and solver status, including provenance for 40 reused fits. Its publication
export revalidates the saved graphs and writes local figures, and complete
operating-point tables. Code, persisted instances, and analysis artifacts will accompany the submission as
anonymized supplementary material.
The code repository is available at 

\url{https://github.com/AtomXT/micodag-intervenion}.

\bibliographystyle{plainnat}
\bibliography{clear2027/references}

\appendix
\section{Expanded Likelihood-Level Reference Formulation}
\label{app:expanded-reference-formulation}

\subsection{Mechanism representation}
\label{app:mechanism-representation}

The mechanism-sharing relation \eqref{eq:mechanism-sharing} can equivalently be expressed
through the structural coefficients and error variances as follows.

Recall that \(\mathcal E\) is the environment set, \(e_0\) is the
observational environment, and
\(\mathcal E_{\setminus0}=\mathcal E\setminus\{e_0\}\). For each
\(e\in\mathcal E\), \(I_e\subseteq[p]\) is its intervention-target set,
with \(I_{e_0}=\varnothing\). Define the diagonal matrix \(R_{I_e}\),
whose \(j\)th diagonal entry is one if node \(j\) is targeted and zero otherwise:
\[
    R_{I_e}
    =
    \operatorname{diag}\bigl(
        \mathbb I[1\in I_e],\ldots,\mathbb I[p\in I_e]
    \bigr).
\]
A stochastic hard intervention replaces the structural equations of the
targeted nodes with exogenous noise and therefore deletes all of their incoming
edges \citep{korb2004varieties}. Therefore, the structural coefficient matrix in
environment \(e\) is \(B^e=B(\mathbf I_p-R_{I_e})\).
For example, if \(I_e=\{k\}\), then the \(k\)th column of \(B\) is set to zero,
removing all incoming edges to node \(k\), while every other column remains
unchanged.  In general, \(I_e\) may contain zero, one, or multiple nodes.

Non-target error variances are shared with the observational mechanism,
whereas a targeted node may have an environment-specific variance. Let
\(\Omega^e\) denote the covariance matrix of the structural errors in
environment \(e\). Its diagonal entries are their variances, and its
off-diagonal entries are zero because the errors are independent. Set
\(\Omega^{e_0}=\Omega\), and, for each \(e\in\mathcal E_{\setminus0}\), define
its positive diagonal entries by
\[
    \Omega^e_{jj}
    =
    \begin{cases}
        \omega_j^e, & j\in I_e,\\
        \Omega_{jj}, & j\notin I_e,
    \end{cases}
    \qquad \omega_j^e>0,
\]
where the separate parameter \(\omega_j^e\) is introduced only when
\(j\in I_e\).
A generic observation in environment \(e\) then satisfies
\[
    X^e=(B^e)^\top X^e+\epsilon^e,
    \qquad
    \epsilon^e\sim\mathcal N(0,\Omega^e),
\]
and has covariance
\(\Sigma^e=(\mathbf I_p-B^e)^{-\top}\Omega^e(\mathbf I_p-B^e)^{-1}\).
Let \(D^e=(\Omega^e)^{-1}\) and define
\begin{equation*}
    \Gamma^e
    =
    (\mathbf I_p-B^e)(D^e)^{1/2}
    =
    \bigl(\mathbf I_p-B(\mathbf I_p-R_{I_e})\bigr)(D^e)^{1/2}.
\end{equation*}
Then \((\Sigma^e)^{-1}=\Gamma^e{\Gamma^e}^{\top}\), and
\[
    \Gamma^e_{\cdot j}
    =
    \begin{cases}
        \Gamma^{e_0}_{\cdot j}, & j\notin I_e,\\
        \sqrt{D^e_{jj}}\,\mathbf e_j, & j\in I_e,
    \end{cases}
\]
where \(\mathbf e_j\) is the \(j\)th standard basis vector.  Thus a non-target
column shares its baseline structural coefficients and diagonal precision
scale, whereas a target column has zero off-diagonal entries and may have an
environment-specific diagonal scale.

Because \(B^e\) is acyclic, \(\det(\mathbf I_p-B^e)=1\). Substituting
\((\Sigma^e)^{-1}=\Gamma^e{\Gamma^e}^{\top}\) into the Gaussian likelihood
gives \(\mathcal L(\Gamma)\)
in Section~\ref{sec:int-formulation}, up to an additive constant and a common
positive factor, with the same empirical second moments and weights.

\subsection{Exact mixed-integer encoding}
\label{app:reference-encoding}

This appendix encodes the model class \(\mathfrak M\) and objective
\eqref{eq:joint-learning-objective} from Section~\ref{sec:int-formulation}.
It is the exact reference problem used by Proposition~\ref{prop:parent-set-equivalence}. The empirical study uses the parent-set formulation rather than this encoding.

Let \(E_{\mathrm{super}}\subseteq\{(i,j):i\neq j\}\) be a directed
superstructure. Use a large enough constant \(M>0\) and penalties
\(\lambda_g,\lambda_\delta\geq0\) from Section~\ref{sec:int-formulation}.  The binary variables \(g_{ij}\) select baseline
edges, the continuous order variables \(\psi_j\) enforce acyclicity, and
\(\gamma_j^e\) denotes the candidate scale
\(\sqrt{D^e_{jj}}=1/\sqrt{\omega_j^e}\) used when node \(j\) is targeted.
Let \(\mathcal F_{\mathrm{graph}}\) collect the superstructure, coefficient
bounds, topological-order constraints, and variable domains listed below.
The exact reference model
is compactly written as
\begin{subequations}
\begin{align}
\min_{\Gamma,\delta,\psi,g,\gamma}
&\mathcal L(\Gamma)
+\lambda_g\sum_{(i,j)\in E_{\mathrm{super}}}g_{ij}
+\lambda_\delta\sum_{e\in\mathcal E_{\setminus0}}\sum_{j=1}^p\delta_j^e
\\
\text{s.t.}\quad
&\Gamma^e_{ij}
=\Gamma^{e_0}_{ij}(1-\delta_j^e)
&&(e\in\mathcal E_{\setminus0};\ i\neq j),
\label{eq:int_target_1}\\
&\Gamma^e_{jj}
=\Gamma^{e_0}_{jj}(1-\delta_j^e)+\gamma_j^e\delta_j^e
&&(e\in\mathcal E_{\setminus0};\ j=1,\ldots,p).
\label{eq:int_target_2}\\
&(\Gamma^{e_0},\delta,\psi,g,\gamma)
\in\mathcal F_{\mathrm{graph}}.
\end{align} 
\label{Problem:interventional}
\end{subequations}
Constraints~\eqref{eq:int_target_1}--\eqref{eq:int_target_2} share the baseline
column when \(\delta_j^e=0\) and replace it by an environment-specific diagonal
scale when \(\delta_j^e=1\).  Their binary--continuous products have exact
McCormick linearizations because all continuous factors are bounded
\citep{mccormick1976computability}; the expansion follows below. Thus
Problem~\eqref{Problem:interventional} is mixed-integer convex representable. 

The condition
\((\Gamma^{e_0},\delta,\psi,g,\gamma)\in
\mathcal F_{\mathrm{graph}}\) in
Problem~\eqref{Problem:interventional} expands to
\begin{align*}
&-Mg_{ij}\leq\Gamma^{e_0}_{ij}\leq Mg_{ij}
&&((i,j)\in E_{\mathrm{super}}),\\
&\Gamma^{e_0}_{ij}=0
&&(i\neq j,\ (i,j)\notin E_{\mathrm{super}}),\\
&\Gamma^{e_0}_{jj}\leq M
&&(j=1,\ldots,p),\\
&1-p+pg_{ij}\leq\psi_j-\psi_i
&&((i,j)\in E_{\mathrm{super}}),\\
&1\leq\psi_j\leq p
&&(j=1,\ldots,p),\\
&g_{ij}\in\{0,1\}
&&((i,j)\in E_{\mathrm{super}}),\\
&\delta_j^e\in\{0,1\}
&&(e\in\mathcal E_{\setminus0};\ j=1,\ldots,p),\\
&\gamma_j^e\leq M
&&(e\in\mathcal E_{\setminus0};\ j=1,\ldots,p).
\end{align*}
The first two families encode the directed superstructure, the next two impose
the common coefficient bounds, $M$, and a topological order, and the remaining
families give the variable domains.  Known or partially known target status is
incorporated by fixing the corresponding \(\delta_j^e\).

For completeness, if \(x\in[\ell,u]\), \(y\in\{0,1\}\), and \(z=xy\), the
exact McCormick linearization is
\[
    \ell y\leq z\leq uy,
    \qquad
    x-u(1-y)\leq z\leq x-\ell(1-y).
\]
For the off-diagonal switching relation
\eqref{eq:int_target_1}, take \(x=\Gamma^{e_0}_{ij}\),
\(y=1-\delta_j^e\), \(z=\Gamma^e_{ij}\), and \(\ell=-M,\ u=M\). For
\eqref{eq:int_target_2}, introduce
\(z_{0j}^e=\Gamma^{e_0}_{jj}(1-\delta_j^e)\) and
\(z_{1j}^e=\gamma_j^e\delta_j^e\), apply the same inequalities with
\(\ell=0,\ u=M\) to each product, and impose
\(\Gamma^e_{jj}=z_{0j}^e+z_{1j}^e\). The diagonal factors retain their
positive domains; the auxiliary products may equal zero when their binary
multiplier is zero. All remaining
constraints are linear, and \(\mathcal L(\Gamma)\) is convex on the positive
diagonal domain, yielding the mixed-integer convex formulation asserted in the
main text.

\subsection{Local-score conditions}
\label{app:local-score-conditions}
The closed-form local scores derived in Section~\ref{sec:local-scores}
are obtained by minimizing over the likelihood parameters without
imposing the bounds in Section~\ref{sec:int-formulation}.
For these scores to equal the minima of the bounded local problems,
we require, for every enumerated triple \((j,A,S)\), that \(Q_{AA}\)
be nonsingular when \(A\neq\varnothing\), that \(r_{jAS}>0\), and
that the fitted baseline parameters satisfy:
\begin{equation*}
    \|\Gamma^{0,*}_{Aj}\|_\infty\leq M,
    \qquad
    0<\Gamma^{0,*}_{jj}\leq M.
\end{equation*}
For an empty parent set, the off-diagonal bound is vacuous.
For every interventional environment \(e\) and node \(j\), we also require
\(\widehat\Sigma^e_{jj}>0\) and
\(1/\sqrt{\widehat\Sigma^e_{jj}}\leq M\).
These conditions ensure that the closed-form fits also solve the bounded local problems. Because the candidate collection is finite, a sufficiently large finite \(M\) contains every local optimum whenever the stated matrix and positivity conditions hold. Increasing \(M\) beyond this threshold leaves the precomputed scores unchanged. Thus, the use of \(M\) establishes equivalence with the bounded reference formulation but does not appear in the final parent-set MIP.

\section{Proof of parent-set equivalence}
\label{app:parent-set-equivalence-proof}

\begin{proof}\textbf{of Proposition~\ref{prop:parent-set-equivalence}.}
Let \(v_{\mathrm{joint}}\) and \(v_{\mathrm{PS}}\) denote the optimal values
of the joint problem~\eqref{eq:joint-learning-objective} and the parent-set
MIP~\eqref{eq:parent-set-mip-objective}--\eqref{eq:parent-set-binary},
respectively. We establish equality by proving both inequalities.

First fix an admissible DAG \(\mathcal G=([p],E)\), and write
\(A_j=\{i:(i,j)\in E\}\) for node $j$'s parent sets. The key identity is
\[
\begin{gathered}
    \min_{\substack{I,\Gamma:\\(\mathcal G,I,\Gamma)\in\mathfrak M}}
    \left\{\mathcal L(\Gamma)+\lambda_g|E|
    +\lambda_\delta\sum_{e\in\mathcal E_{\setminus0}}|I_e|\right\}
    \\
    =\sum_{j=1}^{p}\min_{S\in\mathcal S_j}s(j,A_j,S)
    =\sum_{j=1}^{p}s(j,A_j).
\end{gathered}
\]
For this fixed graph, the likelihood and both penalties separate by node,
with \(|E|=\sum_j|A_j|\). The assumption that feasible target patterns form
\(\prod_j\mathcal S_j\) allows each node's pattern to be chosen independently.
The local-score conditions in Section~\ref{sec:local-scores} ensure that the
closed-form parameter minima satisfy the original bounds and attain the
scores \(s(j,A_j,S)\), establishing the identity.

\emph{From the joint problem to the parent-set MIP.}
Take any feasible solution \((\mathcal G,I,\Gamma)\) of the joint problem
and identify its parent set \(A_j\) for each node. These choices are feasible
for the MIP: the edges belong to \(E_{\mathrm{super}}\), and a topological
ordering of \(\mathcal G\) supplies order variables satisfying the acyclic constraints
\eqref{eq:parent-set-order}--\eqref{eq:parent-set-order-bounds}.
By definition, the total score \(\sum_j s(j,A_j)\) is no
greater than the objective of the original solution of the parent-set formulation, since it minimizes over the remaining target and parameter choices. Hence
\(v_{\mathrm{PS}}\leq v_{\mathrm{joint}}\).

\emph{From the parent-set MIP to the joint problem.}
Take any feasible parent-set solution, with parent sets \(A_j\).
Its acyclic constraints ensure that the edges \(i\to j\) with
\(i\in A_j\) form an admissible DAG. For each node, we have the saved
minimizing pattern \(S^*(j,A_j)\) and its fitted parameters from
Section~\ref{sec:local-scores}, then we can reconstruct the parameter in the orignal joint problem:

For the selected parent sets \(A_j\), set
\(I_e=\{j:e\in S^*(j,A_j)\}\), recover the optimized columns, and copy each
baseline column into its non-target environments. Recover the original
parameters using
\[
    B_{ij}=-\Gamma^0_{ij}/\Gamma^0_{jj}\quad(i\neq j),\qquad
    \Omega_{jj}=1/(\Gamma^0_{jj})^2,\qquad
    \omega_j^e=1/(\Gamma^{e,*}_{jj})^2\quad(j\in I_e).
\]
These choices combine into a
feasible joint solution: the target patterns are jointly feasible, the
parameters satisfy the original bounds, and the columns obey the
mechanism-sharing rule. Its objective is exactly \(\sum_j s(j,A_j)\),
the objective of the parent-set solution. Hence
\(v_{\mathrm{joint}}\leq v_{\mathrm{PS}}\).

Together, the inequalities give \(v_{\mathrm{PS}}=v_{\mathrm{joint}}\).
Applying the second construction to an optimal parent-set solution therefore
recovers an optimal graph, target assignment, and likelihood parameters for
the joint problem. The same conclusion holds for its reference
encoding~\eqref{Problem:interventional}, whose equivalence to
\eqref{eq:joint-learning-objective} is established in
Appendix~\ref{app:reference-encoding}.
\end{proof}

\section{Algorithm}
\label{app:parent-set-implementation}

Algorithm~\ref{alg:ps-mip} precomputes optimized local scores, selects parent
sets subject to acyclicity, and reconstructs the model. It assumes the
local-score conditions in Section~\ref{sec:local-scores}.

\begin{algorithm2e}[tbp]
\DontPrintSemicolon
\caption{PS-MIP: joint estimation through optimized local scores}
\label{alg:ps-mip}
\KwIn{\(\{(\widehat\Sigma^e,w_e)\}_{e\in\mathcal E}\), observational
environment \(e_0\), candidate edges \(E_{\mathrm{super}}\), penalties
\((\lambda_g,\lambda_\delta)\), parameter bounds, and compatible patterns
\(\{\mathcal S_j\}_{j=1}^p\).}
\KwOut{Graph \(\widehat{\mathcal G}\), targets \(\widehat I\), fitted
parameters, and optimization diagnostics.}
\BlankLine
\textbf{1. Local-score preprocessing}\;
\For{\(j\in[p]\), \(A\subseteq\operatorname{pa}^{\mathrm{sup}}(j)\)}{
    \For{\(S\in\mathcal S_j\)}{
        Compute \(s(j,A,S)\) and its fitted parameters using
        Section~\ref{sec:local-scores}\;
    }
    Choose \(S^*(j,A)\in\arg\min_{S\in\mathcal S_j}s(j,A,S)\)\;
    \(s(j,A)\leftarrow s\bigl(j,A,S^*(j,A)\bigr)\)\;
    Store \(S^*(j,A)\) and its fitted parameters\;
}
\BlankLine
\textbf{2. Global parent-set selection}\;
Solve \eqref{eq:parent-set-mip-objective}--\eqref{eq:parent-set-binary}
with scores \(s(j,A)\), obtaining a feasible solution \(\widehat x\)\;
\(\widehat A_j\leftarrow\text{the unique }A\text{ with }\widehat x_{jA}=1\),
\(j\in[p]\)\;
\BlankLine
\textbf{3. Reconstruction}\;
\(\widehat E\leftarrow\{(i,j):i\in\widehat A_j\}\),
\(\widehat{\mathcal G}\leftarrow([p],\widehat E)\)\;
\(\widehat I_e\leftarrow\{j:e\in S^*(j,\widehat A_j)\}\),
\(e\in\mathcal E\)\;
Recover coefficients and variances from the stored parameters using
the reconstruction paragraph in Appendix~\ref{app:parent-set-equivalence-proof}\;
\end{algorithm2e}

\section{Additional Screening and Experiment Details}
\label{app:synthetic-diagnostics}

\subsection{Observational screening justification}
\label{app:screening-justification}
In this section, we justify the use of the observational precision matrix as a conservative superstructure screen; the specific graphical-lasso implementation and parameter choices used in the experiments are given in Appendix~\ref{app:experimental-protocol}.

Let \(\Theta^0=(\mathbf I_p-B)D(\mathbf I_p-B)^\top\) be the observational precision matrix,
where \(D=\Omega^{-1}\) is diagonal.  For \(i\neq j\),
\[
    \Theta^0_{ij}
    =
    -D_{jj}B_{ij}-D_{ii}B_{ji}
    +\sum_{k=1}^p D_{kk}B_{ik}B_{jk}.
\]
Thus \(\Theta^0_{ij}=0\) whenever \(i\) and \(j\) are neither adjacent in the
DAG nor parents of a common child, so the support graph of \(\Theta^0\) is
contained in the moral graph.  Under the no-cancellation condition of
\citet{loh2014high}, the edge sets coincide; this condition holds generically
when the nonzero coefficients vary continuously.

Let $\widehat\Theta^0$ be an estimate of the observational precision
matrix, and retain an undirected pair $\{i,j\}$ whenever
$\widehat\Theta^0_{ij}\neq 0$. Suppose
\[
    \|\widehat\Theta^0-\Theta^0\|_{\max}\leq \epsilon
\]
and
\[
    \min_{\{i,j\}\in E(M(\mathcal G))}
    |\Theta^0_{ij}|>\epsilon, \text{ where } M(\mathcal G) \text{ the moral graph of a DAG } \mathcal G.
\]
Then every edge of the moral graph is retained by the screen. Indeed, for
any $\{i,j\}\in E(M(\mathcal G))$,
\[
    |\widehat\Theta^0_{ij}|
    \geq
    |\Theta^0_{ij}|
    -
    |\widehat\Theta^0_{ij}-\Theta^0_{ij}|
    >0.
\]
Hence the retained graph contains the moral graph and, in particular,
contains the skeleton of the true DAG. Including both orientations of each
retained pair therefore gives a candidate directed superstructure containing the true DAG. 

Elementwise graphical-lasso results provide this type of high-probability
screening guarantee under standard sparsity, signal-strength, tail, and
incoherence conditions \citep{ravikumar2011high}.  The experiment uses the
estimate only as a conservative one-sided screen.

\subsection{Experimental protocol}
\label{app:experimental-protocol}

All experiments are implemented in Python 3.9 and solved using Gurobi 12.0.3 through its Python API. We run the computations on a single compute node with an Intel Xeon Platinum 8592+ processor, using 8 CPU cores and 16 GB of RAM. Unless otherwise specified, Gurobi parameters are kept at their default values. Time limits are set to 3600 seconds.  PS-MIP screening and local-score
precomputation occur before and outside the Gurobi limit, but are included
in reported fitting time.

\paragraph{Synthetic data generation.}
For each combination of \(p\in\{10,20,30\}\) and edge-density parameter
\(e\in\{1,2\}\), defined by \(\mathbb E[|E|]=ep\), we generate ten
independent ordered random DAGs.
Each forward edge is included independently with probability
\(\min\{1,2e/(p-1)\}\), giving \(ep\) expected edges.
Nonzero coefficients are sampled from \(\{-0.8,-0.6,0.6,0.8\}\) and baseline
noise variances from \(\{1,2,4\}\). Each instance has 1000 observational
samples and five interventional environments with 200 samples each.
In each intervention, the target count is uniform on
\(\{1,\ldots,\lfloor p/2\rfloor\}\); targets are chosen without replacement,
their incoming edges removed, and their noise variances sampled uniformly
from \([1,4]\). Graphs, coefficients, targets, and data are independently
generated across replicates: 60 saved instances with 2000 samples each.

\paragraph{Screening and parameter settings.} Following the screening construction justified in \linebreak[4]Appendix~\ref{app:screening-justification}, PS-MIP uses observational samples to estimates the observational precision matrix using graphical lasso with
\(\alpha_{\mathrm{screen}}=\sqrt{\log(p)/n_{\mathrm{obs}}}\), where $n_{\mathrm{obs}}$ is the number of observational samples.  We retain every pair $\{i,j\}$ for which the estimated off-diagonal
precision entry is nonzero, and include both orientations of each
retained pair in the candidate superstructure.
The screen retains 91.9--99.4\% of true directed edges, pooled within each
combination of graph size and density. 

\ignore{
\begin{table}[H]
    \caption{PS-MIP screening diagnostics.  True-edge recall is pooled
    across directed edges in the persisted replicates; screen edges and
    candidate parent sets are replicate means.}
    \label{tab:ps-mip-screen-diagnostics}
    \centering
    \begingroup
    \small
    \setlength{\tabcolsep}{5pt}
    \begin{tabular}{@{}rrrrr@{}}
        \toprule
        $p$ & $e$ & $n$ & True-edge recall & Mean screen edges / parent sets \\
        \midrule
        10 & 1 & 10 & 0.970 & 22.1 / 440.5 \\
        10 & 2 & 10 & 0.919 & 31.9 / 1141.6 \\
        20 & 1 & 10 & 0.989 & 50.5 / 3834.5 \\
        20 & 2 & 10 & 0.950 & 90.3 / 40581.5 \\
        30 & 1 & 10 & 0.994 & 87.8 / 18943.1 \\
        30 & 2 & 10 & 0.956 & 167.2 / 990604.7 \\
        \bottomrule
    \end{tabular}
    \par\smallskip
    \begin{minipage}{\linewidth}
    \footnotesize The screen is recomputed from the persisted observational
    sample using $\alpha=\sqrt{\log(p)/1000}$ and the experiment's exact
    GraphicalLasso support rule.
    \end{minipage}
    \endgroup
\end{table}
}

With \(n=2000\), unknown-target PS-MIP uses equal graph and target
penalties from \(\{0.5,1,2\}\allowbreak\log(n)/n\), and oracle PS-MIP uses the same graph
penalties with fixed target penalties since the targets are known.  GnIES and GIES use
\(\{0.25,0.5,1,2,4\}\) times their native \(\tfrac12\log(n)\) Bayesian Information Criterion (BIC) penalty.
For UT-IGSP and IGSP, the conditional-independence and invariance levels are
tied over \(\{10^{-5},\allowbreak10^{-4},\allowbreak10^{-3},\allowbreak10^{-2},\allowbreak0.05\}\); the search depth is 4,
and UT-IGSP uses 10 randomized initial runs.  BaCaDI uses one predeclared
configuration per instance with 20 particles and 3000 Stein Variational Gradient Descent (SVGD) steps.  DCDI-G is
attempted separately under the same fit budget.

\paragraph{Graph metrics.}
\label{app:graph-metrics}
We convert fitted graphs and estimated targets to interventional CPDAGs.
Let \(\widehat{\mathcal C}\) and \(\mathcal C^\star\) contain the DAG
extensions of the estimated and true classes, and \(E(D)\) denote a DAG's
directed edges. Following \citet{taeb2024learning}, we use
\begin{align*}
\operatorname{FDP}
&=
\max_{\widehat D\in\widehat{\mathcal C}}
\min_{D^\star\in\mathcal C^\star}
\frac{\left|E(\widehat D)\setminus E(D^\star)\right|}
     {\left|E(\widehat D)\right|},
\\
\operatorname{TDP}
&=
\min_{D^\star\in\mathcal C^\star}
\max_{\widehat D\in\widehat{\mathcal C}}
\frac{\left|E(D^\star)\cap E(\widehat D)\right|}
     {\left|E(D^\star)\right|}.
\end{align*}
FDP measures the worst-case fraction of estimated edges unmatched by any
true-class DAG; TDP measures the worst-case fraction of true edges recovered
by the best matching estimated-class DAG. An empty reference edge set has
missing-edge proportion zero: FDP is \(0\) for an empty estimated graph,
and TDP is \(1\) for an empty true graph.
We report the average of metrics across all replicates.

\paragraph{Target metrics.}
For replicate \(r\), let
\(\Delta^{(r)},\widehat\Delta^{(r)}\in\{0,1\}^{5\times p}\) be the generating
and estimated environment-by-node target indicators.  The standard target
true- and false-positive rates are
\begin{align*}
\operatorname{TPR}^{(r)}
&=
\frac{\sum_{k=1}^{5}\sum_{j=1}^{p}
      \mathbb I\{\widehat\Delta_{kj}^{(r)}=1,
                    \Delta_{kj}^{(r)}=1\}}
     {\sum_{k=1}^{5}\sum_{j=1}^{p}
      \mathbb I\{\Delta_{kj}^{(r)}=1\}},
\\
\operatorname{FPR}^{(r)}
&=
\frac{\sum_{k=1}^{5}\sum_{j=1}^{p}
      \mathbb I\{\widehat\Delta_{kj}^{(r)}=1,
                    \Delta_{kj}^{(r)}=0\}}
     {\sum_{k=1}^{5}\sum_{j=1}^{p}
      \mathbb I\{\Delta_{kj}^{(r)}=0\}}.
\end{align*}

\subsection{Computational Time}

\begin{table}[t]
\centering
\small
\caption{Median fitting time (seconds) across replicates and tuning
settings, grouped by graph size \(p\) and density parameter \(e\). TL denotes a median fitting time hit
one-hour limit; a dash indicates no qualifying result.}
\label{tab:synthetic-runtime}
\setlength{\tabcolsep}{5pt}
\begin{tabular}{@{}lrrrrrr@{}}
\toprule
& \multicolumn{2}{c}{\(p=10\)}
& \multicolumn{2}{c}{\(p=20\)}
& \multicolumn{2}{c}{\(p=30\)} \\
\cmidrule(lr){2-3}
\cmidrule(lr){4-5}
\cmidrule(lr){6-7}
Method
& \(e=1\) & \(e=2\)
& \(e=1\) & \(e=2\)
& \(e=1\) & \(e=2\) \\
\midrule
PS-MIP
& 1.43 & 7.74 & 13.73 & TL & 57.02 & TL \\
UT-IGSP
& 0.15 & 0.32 & 0.27 & 0.48 & 0.57 & 0.67 \\
GnIES
& 22.33 & 66.29 & 739.98 & 2285.30 & 3234.33 & -- \\
\midrule
PS-MIP (oracle)
& 0.28 & 0.94 & 1.43 & 745.09 & 7.57 & TL \\
IGSP (oracle)
& 0.08 & 0.22 & 0.34 & 2.07 & 1.62 & 15.10 \\
GIES (oracle)
& 0.34 & 0.51 & 3.74 & 7.45 & 18.73 & 40.83 \\
\bottomrule
\end{tabular}
\end{table}

Table~\ref{tab:synthetic-runtime} reports median fitting times across
graph sizes and densities. PS-MIP has median fitting times below one
minute in all \(e=1\) settings, but its computational cost increases
substantially for denser graphs, exceeding one hour for \(e=2\) with
\(p\in\{20,30\}\). Providing the true intervention targets reduces
PS-MIP's median fitting time in every setting, although the largest
dense setting still exceeds one hour. UT-IGSP, oracle IGSP, and oracle
GIES are substantially faster than PS-MIP in the larger dense settings.
These results indicate that computational cost remains a limitation
of PS-MIP as graph size and density increase. When a fit reaches the time limit, PS-MIP returns the best incumbent found so far---a feasible DAG that is not certified optimal---rather than failing; such settings are reported as time-limited (TL) and marked with crosses in the corresponding figures.

\paragraph{Comparison with the original joint formulation.}
We compare PS-MIP with the original joint MIP
\eqref{Problem:interventional}, implemented with the exact switching-variable
linearization in Appendix~\ref{app:reference-encoding}. Both formulations use
the same 60 synthetic instances, three penalty settings, unknown intervention
targets, observational screening rule, and parameter bounds. Each fit uses
eight solver threads, a one-hour optimization limit, and a relative MIP-gap
tolerance of \(10^{-4}\). Table~\ref{tab:joint-mip-comparison} reports
median total fitting times, including preprocessing and model construction,
together with optimal statuses and relative gaps as reported by Gurobi.

\begin{table}[t]
\centering
\small
\caption{PS-MIP and the original joint MIP with unknown intervention targets.
Time is median total fitting time in seconds; preparation is included and lies
outside the one-hour solver limit. Optimal statuses and relative gaps are
Gurobi-reported; gaps are median percentages over available feasible fits.
Each \((p,e)\) cell has 30 planned fits (ten replicates and three penalties).
PS-MIP has 26 available fits at \(p=30,e=2\) and 30 in every other cell;
the joint MIP has 30 in every cell. Reported-optimal counts use the planned
30 fits as their denominator. The relative gap is the solver's certificate of suboptimality: the percentage difference between the objective of the returned solution and a proven lower bound on the best achievable objective. A $0\%$ gap certifies a globally optimal solution; a gap of, e.g., $9.78\%$ guarantees the returned objective is within $9.78\%$ of the global optimum. \emph{Reported opt.}\ counts the fits the solver certified optimal (gap within the $10^{-4}$ tolerance); a solver may return a feasible solution without certifying its optimality. The gap measures objective function suboptimality, not causal-recovery error.}
\label{tab:joint-mip-comparison}
\setlength{\tabcolsep}{3pt}
\begin{tabular}{@{}rr rrr rrr@{}}
\toprule
 & & \multicolumn{3}{c}{PS-MIP} & \multicolumn{3}{c}{Original joint MIP} \\
\cmidrule(lr){3-5}\cmidrule(lr){6-8}
\(p\) & \(e\) & Time (s) & Reported opt. & Gap (\%)
& Time (s) & Reported opt. & Gap (\%) \\
\midrule
10 & 1 & 1.43 & 30/30 & 0.00 & 577.22 & 21/30 & 0.01 \\
10 & 2 & 7.74 & 30/30 & 0.00 & TL& 5/30 & 2.22 \\
20 & 1 & 13.73 & 30/30 & 0.00 & TL& 2/30 & 2.84 \\
20 & 2 & TL & 10/30 & 1.30 & TL& 2/30 & 16.61 \\
30 & 1 & 57.02 & 26/30 & 0.00 & TL& 4/30 & 9.78 \\
30 & 2 & TL & 0/30 & 5.61 & TL& 0/30 & 23.65 \\
\bottomrule
\end{tabular}
\end{table}

\section{Causal Chambers Protocol and Additional Results}
\label{app:chambers}

\subsection{Data and experimental protocol}
\label{app:chamber-protocol}

\paragraph{Data and reference graph.}
We use configuration \texttt{scm\_4} of the light-tunnel dataset
\linebreak[4]\texttt{lt\_camera\_v1} \citep{gamella2024causal}, with environment order
reference, red, green, blue, and 10,000 observations per environment.
The variable order is red, green, blue, current, ir\_1, vis\_1, ir\_2,
vis\_2, ir\_3, vis\_3. Move details are provided in
the accompanying documentation.\footnote{\raggedright
Dataset and protocols: \url{https://github.com/juangamella/causal-chamber/tree/main/datasets/lt_camera_v1}.\\
Reference adjacency: \url{https://github.com/juangamella/causal-chamber-package/tree/main/causalchamber/ground_truth/adjacencies}.}

The official camera adjacency restricted to these variables supplies seven
outgoing physical links from each of \(R,G,B\), for 21 physical edges.
The programmed actuator chain adds \(R\to G\) and \(G\to B\).
Omitted documented parents, including polarizer angles, auxiliary light
settings, and sensor parameters, are verified to be fixed across all four
environments. The protocol samples truncated-normal exogenous actuator
values, propagates them through its programmed equations, and sends
integer-valued intensity commands to the physical chamber. An intervention
zeros incoming coefficients of its chosen actuator and shifts its exogenous
distribution; the seven sensor channels remain real measurements. This
protocol specifies how the interventions remove incoming effects at the
targeted actuators.

\paragraph{Preprocessing and target information.}
For each variable \(j\), the reference mean \(\bar x_{0j}\) and population
standard deviation \(s_{0j}\) define
\[
    z_{eij}=(x_{eij}-\bar x_{0j})/s_{0j}
\]
in every environment \(e\). All methods
receive these same arrays. PS-MIP centers each environment; GnIES and GIES
use their documented internal environment centering. UT-IGSP and IGSP retain
the arrays' mean shifts for invariance testing.

PS-MIP receives no fixed interventional target labels; UT-IGSP receives no
known targets; GnIES receives neither the set of variables targeted across
environments nor an initial guess for that set. All methods receive the
identity of the observational environment. Oracle GIES and IGSP additionally
receive the documented target lists \(\{R\},\{G\},\{B\}\).

\paragraph{Prespecified settings and computation.}
We fixed the parameter settings before fitting. With \(n=40{,}000\) total observations, PS-MIP uses the 17 graph penalties
\[
    \lambda_g=2^{k/2}\log(n)/n,\qquad k=-8,\ldots,8,
    \qquad \lambda_\delta=16\log(n)/n.
\]
The target penalty is fixed across environments and graph-penalty settings.
No screening is
used: all 45 pairs and 5,120 parent sets remain available. GnIES and GIES use
15 multipliers \(2^k\), \(k=-4,\ldots,10\), of their native
\(\tfrac12\log n\) penalty. \linebreak[4]UT-IGSP and IGSP use Gaussian conditional-independence
(CI) and invariance tests, CI significance levels \(0.001,0.01,0.05,0.1,0.2,0.3,0.4,0.5\), fixed invariance
level \(10^{-5}\), depth four, and ten starts. The solver/search seed is
20260901, with eight solver threads and one numerical-library thread.
The Quest environment uses Python 3.9.5, Gurobi 12.0.3, GnIES 0.3.3,
GIES 0.0.3, NumPy 1.26.3, and SciPy 1.11.4.

Each PS-MIP solve has a 3,600-second solver limit; score precomputation
occurs outside that limit but it is included in the reported elapsed times
(Table~\ref{tab:chamber-completion}). Baselines have a 3,600-second limit
for their fitting procedure. 

\subsection{Evaluation and additional results}
\label{app:chamber-results}

\begin{figure}[t]
    \centering
    \includegraphics[width=0.75\linewidth]{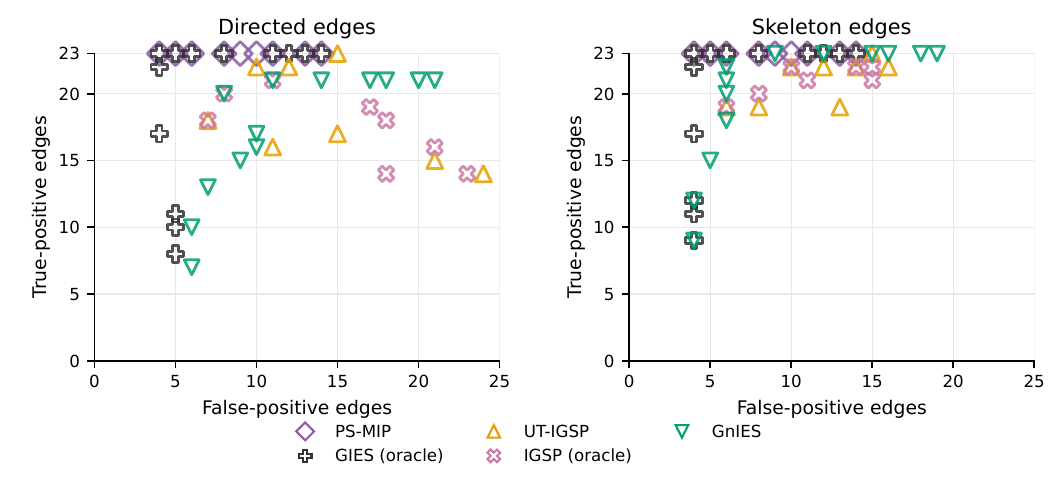}
    \caption{Five-method Chamber comparison. Targets are withheld from
    PS-MIP, UT-IGSP and GnIES, and supplied as documented oracle lists to
    GIES and IGSP. The false-positive axis spans a fixed 0--25 range in every panel. All 63 fits are represented: exact duplicates within a
    method and panel are drawn once, producing 51 directed and 50 skeleton
    markers, with coincident points from different methods kept separate.
    All points fall within this range.
    Table~\ref{tab:chamber-completion} reports solver status.}
    \label{fig:chamber-five}
\end{figure}

Figure~\ref{fig:chamber-five} adds oracle GIES and IGSP to the unknown-target
comparison in the main text. PS-MIP recovers all 23 reference directions
across its tested path, with four to fourteen extra edges. Oracle GIES
matches its sparsest point: both methods attain 23 TP/4 FP for directed
and skeleton recovery. At these points, all 23 true directions are shared
by every DAG in the saved interventional equivalence classes; two of the
four extra adjacencies have ambiguous orientations.

Table~\ref{tab:chamber-completion} reports all 63 planned fits, each of which
returned a validated graph. All 17 PS-MIP runs reached the time limit with
feasible solutions. Their relative MIP
gaps range from 10.81\% to 16.57\%, with a median of 13.78\%. These gaps quantify remaining
objective uncertainty using the solver's bound, rather than causal recovery
error.

\begin{table}[t]
\centering
\small
\caption{Chamber completion and computational diagnostics. All planned fits
returned an estimate; no failures or missing results are omitted. Times cover
data loading, method execution/precomputation, and result validation within each per-setting run, excluding run startup and final file writing. The total sums across
settings; reused fits retain their original times. Returned gives completed
runs out of planned runs. Nonoptimal counts returned PS-MIP solutions without
certified optimality. Median (s) is the median elapsed time per run in seconds. A dash
means that a comparable optimality certificate is unavailable, not a zero gap.}
\label{tab:chamber-completion}
\begin{tabular}{lrrrrr}
\toprule
Method & Returned & Nonoptimal & Median (s) & Total (min) & Median gap (\%) \\
\midrule
PS-MIP & 17/17 & 17 & 3617.71 & 1025.773 & 13.78 \\
GIES (oracle) & 15/15 & -- & 1.00 & 0.235 & -- \\
UT-IGSP & 8/8 & -- & 0.69 & 0.112 & -- \\
IGSP (oracle) & 8/8 & -- & 0.19 & 0.027 & -- \\
GnIES & 15/15 & -- & 100.18 & 24.585 & -- \\
\bottomrule
\end{tabular}
\end{table}

\end{document}